\documentclass [12pt]{article}
\usepackage[dvips]{graphics}
\usepackage{amsmath}
\usepackage{amssymb}
\usepackage{amscd}
\usepackage{amsthm}
\usepackage{amsopn}
\usepackage{xspace}
\usepackage{verbatim}
\usepackage{amsmath}
\usepackage{amsfonts}
\usepackage{colortbl}
\usepackage{hyperref}
\usepackage{amsmath}
\usepackage{color}
\usepackage{srcltx}
\usepackage[utf8]{inputenc}
\newtheorem{theorem}{Theorem}
\newtheorem{proposition}{Proposition}[section]

\newtheorem{remark}{Remark}[section]
\newtheorem{lemma}{Lemma}[section]

\newcommand\AMSname{AMS subject classifications}

\renewcommand{\ln}{\log}

\usepackage{soul}

\let\div\relax 
\DeclareMathOperator{\div}{div}
\DeclareMathOperator{\diam}{diam}

\begin{document}




 



\title{A system with weights and with critical Sobolev exponent}

\author{Asma Benhamida\footnote{University of Tunis El Manar, Faculty of sciences of Tunis, Department of Mathematics, Tunisia.
e-mail : asma.benhamida2019@gmail.com} \, and $\,$ Rejeb Hadiji\footnote{Corresponding author : Université Paris-Est, LAMA (UMR 8050), UPEM, UPEC, CNRS, Cr\'eteil, France. e-mail : rejeb.hadiji@u-pec.fr } }
\date{}

\maketitle

\vspace{0.5cm}

\begin{abstract}
In this paper, we investigate the minimization problem :
$$ \inf_{ \displaystyle{\begin{array}{lll} u \in  H_0^1(\Omega), v \in  H_0^1(\Omega),\\ \quad \| u \|_{L^{q}} =1, \quad \| v \|_{L^{q}} = 1 \end{array}}} \left[ \frac{1}{2} \int_{\Omega} a(x) \vert \nabla u(x) \vert^2dx  + \displaystyle{ \frac{1}{2} \int_{\Omega} b(x) \vert \nabla v (x)\vert^2dx } -  \lambda \displaystyle{\int_{\Omega} u(x)v (x)dx} \right]   $$
where $q=\frac{2N}{N-2}$, $ N \geq 4$, $a$ and $b$ are two continuous positive weight functions. We show the existence of solutions of the previous minimizing problem under some conditions on $a$, $b$, the dimension of the space and the parameter $\lambda$.

\bigskip

\noindent2010 \AMSname : 35A01, 35A15, 35J57, 35J62.
 \medskip

\medskip

\end{abstract}

 \section{Introduction}
 Let $ \Omega \subset \mathbb{R}^N $ be a given smooth bounded domain, with $N \geq 4$. Throughout this paper, we are concerned with the following nonlinear minimization problem :
 \begin{eqnarray}  \label{1.1}
   Q_\lambda = \inf_{ \displaystyle{ \begin{array}{lll} (u, v) \in \left( H_0^1(\Omega) \right)^2 \setminus \{(0,0) \} \end{array}}}  E_{\lambda}(u, v),
 \end{eqnarray}
with
$$ E_{\lambda}(u, v)= \displaystyle {\frac{ 1} {2 \| u \|_{q}^2}{ \int_{\Omega} a(x) \vert \nabla u (x)\vert^2dx}} + \displaystyle{ \frac{ 1}{2 \| v(x) \|_{q}^2} \int_{\Omega} b(x) \vert \nabla v(x) \vert^2dx} - \displaystyle{ \frac{\lambda}{ \| u \|_{q} \| v \|_{q} } \int_{\Omega} u(x)v(x) dx }.  $$
where $a$ and $b$ are positive continuous functions on $\overline{\Omega}$, $\lambda$ is a real constant and $ q=\frac{2N}{N-2}$ is the critical exponent for the Sobolev embedding
\begin{equation} \label{AA}
H_0^1(\Omega) \hookrightarrow L^{q}(\Omega).
\end{equation}

Note that  positive minimizers $u$ and $v$  for \eqref{1.1} are non trivial solutions of the  Euler-Lagrange equation associated to \eqref{1.1} namely :
\begin{equation} \label{eqation-Euler}
 \left\{
       \begin{array}{rclll}
 - \div(a(x) \nabla u) - \lambda v  &=& \Lambda_1 u^{2^{*}-1}  &\mbox{in} & \Omega \\
 - \div(b(x) \nabla v) - \lambda u  &= & \Lambda_2  v^{2^{*}-1} &\mbox{in} &\Omega \\
 u \geq 0, \quad   v &\geq& 0 &\mbox{ in } &\Omega \\
  u =v&=& 0  &\mbox{on}  &\partial \Omega,\\
 \end{array}
    \right.
\end{equation}
 where  $\Lambda_1, \Lambda_2 \in \mathbb{R}$ are the Lagrange multipliers associated to  \eqref{1.1} .


 In the well-known article \cite{BN}, Brezis and Nirenberg treated the problem ~\eqref{1.2} in the special case where the weights $a$ and $b$ are positive constant function. They proved that the problem  has at least one positive solution for  $0 < \lambda< \lambda_1$ when $N\geq 4$ and for $\lambda^{*}< \lambda< \lambda_1$ when $N=3, $ where $\lambda_1$ denotes  the first eigenvalue of $- \Delta$ with homogeneous Dirichlet conditions  and $\lambda^{*}$ is a positive constant. \\
 In the presence of a  non-constant positive and bounded weight $h$, the scalar problem  has been addressed  in  \cite{HY}. See also \cite{HMPY} where the authors consider the existence of minimizers solutions  for the equation
\begin{equation} \label{1.2}
- \div(h(x) \nabla u) = \lambda u + \vert u \vert^{q- 2} u  \quad \mbox{ in }  \Omega,  \quad  u > 0 \mbox{ in }  \Omega, \quad u = 0 \quad \partial\Omega.
\end{equation}
They proved that the existence of a solution depends on the first eigenvalue $ \lambda_{1}^{ h} $ of $ - \div(h(x) \nabla \cdot)$ in $H_0^1(\Omega)$, on the behavior of the function $h$ in the vicinity of its minima, and on the geometry of the domain $\Omega$  In \cite{FS}, Furtado and Souza generalized the problem~\eqref{1.2},  considering a non-homogeneous term $ \vert u \vert^{q-2} u. $ We refer also to \cite{BHVY, H, HV} for more general weights which depend on $u$ and $x$. There are other results in the scalar case for this kind of problem, see the references \cite{H1, H, HV}.  \\
In \cite{AMS}, the authors consider a critical and subcritical system with critical non-linearity term $u^{\alpha-1} v^\beta$, where $\alpha > 1$ and $\beta > 1$. \\
The main difficulty faced when dealing with this problem is the lack of compactness of the embedding~\eqref{AA}. We overcome this problem using the presence of  linear perturbation  terms $\lambda u$ and $\lambda v$. In \cite{CP}, the authors considered an arbitrary  system with $l$ components, $l\geq 2$, for the Yamabe equation on a closed Riemannian manifold; see also~\cite{CTV}.

Note that the shape of the domain can have a strong influence on the type of results one can expect. See for example the seminal work of J.M.~Coron \cite{C}.

Geometric  and physical motivations, in particular in relation to the Yamabe problem, can be found for example in \cite{A}, \cite{CNV} \cite{J}, \cite{K},  \cite{L},  \cite{LP}.

In \cite{K}, the authors connect the Ginzburg-Landau functional used in the physics literature with the Yamabe functional known in mathematics. It is shown  that, if properly interpreted, both functionals upon minimization produce the same
Ginzburg-Landau-type equations used on critical phenomena.

 In this paper, we consider the case where $a$ and $b$ are non constant distinct weights. We prove the  existence of positive solutions $(u, v) $ which depend, among others, on the behaviors of the weights $a(\cdot)$ and $b(\cdot)$ near their minima and the dimension of the space.

 In order  to  state  the problem  and to announce our main results, we introduce some preliminaries. \\
 Let us assume the existence of $x_0, $ in $ \Omega $ such that, in a neighborhood of $x_0, $ $a$ and $b$ behave like
  \begin{equation}  \label{1.3}
     a(x) = a(x_0) + {A_k \vert x-x_0 \vert}^k + {\vert x-x_0 \vert}^k \theta_{a}(x), \quad \mbox{as} \  x \rightarrow x_0,
     \end{equation}
 \begin{equation}  \label{1.4}
 b(x) = b(x_0) + {B_l \vert x-x_0 \vert}^l + {\vert x-x_0 \vert}^l \theta_{b}(x), \quad \mbox{as} \  x \rightarrow x_0,
 \end{equation}
 where $k > 0 $, $ l > 0 $ and $A_k, $ $ B_l $ are positive constants, $\theta_{a}(x) $ and $\theta_{b}(x)$ tend to $0$ when x tends to $x_0$.

  The parameters $k$ and $l$ will play a critical role in the study of our problem. Indeed, if $N \geq 4$ the case $(k, l )$ with $k > 2$ and $l > 2$ is treated through a classical procedure.

 If $0<k\leq 2$, $0<l\leq 2$  the problem is more delicate, we have to assume that the functions $a$ and $ b $ satisfy the following additional conditions :
 \begin{equation}   \label{1.5}
 k A_k \leq \frac{\tilde{a}(x)}{{\vert x - x_0 \vert}^k} \quad \mbox{a.e.} \quad x \in \Omega,
\end{equation}
and
  \begin{equation}   \label{1.6}
 l B_l \leq \frac{\tilde{b}(x)}{{\vert x - x_0 \vert}^l} \quad \mbox{a.e.} \quad x \in \Omega,
\end{equation}
where
$$ 
\tilde{a}(x):= \nabla a(x) \cdot (x-x_0) \quad \mbox{and} \quad  \tilde{b}(x):= \nabla b(x) \cdot (x-x_0).
$$

In order to highlight the difficulty, we consider the blow-up of $ u, v \in H_0^1(\Omega) $ around $ x=x_0$, see \cite{H, HV}. Depending on the parameters $k$, $l$ and $\lambda$ different situations occur in the blow-up scale around the point where the weights reach their minimum. 
More precisely, we consider the function $w_{\varepsilon} $ and $ z_{\varepsilon} $ defined by :
\begin{eqnarray}
 \forall \ \varepsilon > 0, &\ u(x) =  \varepsilon^{-\frac{(N-2)}{2}} w_{\varepsilon} \left( \frac{x-x_0}{\varepsilon} \right), \nonumber \\
 &v(x)= \varepsilon^{-\frac{(N-2)}{2}} z_{\varepsilon} \left( \frac{x-x_0}{\varepsilon} \right).
\end{eqnarray}
One has $w_{\varepsilon}, z_{\varepsilon} \in H_0^1(\Omega_{\varepsilon}) $ with $\Omega_{\varepsilon} = \left\lbrace \varepsilon^{-1} y, \ y \in \Omega \right\rbrace $ and $  \| w_{\varepsilon} \|_{L^{q}( \Omega_{\varepsilon})} = \| z_{\varepsilon} \|_{L^{q}( \Omega_{\varepsilon})}= \| u \|_{L^{q}(\Omega)} = \| v \|_{L^{q}(\Omega)} $, thus

$$E_\lambda(u, v) = \displaystyle{\frac{1}{2} \int_{\Omega_{\varepsilon}} a(\varepsilon y+x_0) \vert \nabla w_{\varepsilon}(y) \vert^2dy} + \displaystyle{\frac{1}{2} \int_{\Omega_{\varepsilon}} b(\varepsilon y+x_0) \vert \nabla z_{\varepsilon}(y) \vert^2dy} -  \lambda \varepsilon^{2} \displaystyle{\int_{\Omega_{\varepsilon}} w_{\varepsilon}(y) z_{\varepsilon}(y) dy}. $$
Consequently, the blow-up around $ x=x_0$ gives
$$ E_\lambda(u, v) \underset{\varepsilon \rightarrow 0}{\sim}   \displaystyle{\frac{1}{2} \int_{\Omega_{\varepsilon}} a(x_0) \vert \nabla w_{\varepsilon}(y) \vert^2dy} + \displaystyle{\frac{1}{2} \int_{\Omega_{\varepsilon}} b(x_0) \vert \nabla z_{\varepsilon}(y) \vert^2dy} -  \lambda \varepsilon^{2} \displaystyle{\int_{\Omega_{\varepsilon}} w_{\varepsilon}(y) z_{\varepsilon}(y) dy} +$$
 $$ \displaystyle{\frac{A_k}{2} \int_{\Omega_{\varepsilon}}\vert \varepsilon y+x_0\vert^k \vert \nabla w_{\varepsilon}(y) \vert^2dy} + \displaystyle{\frac{B_l}{2}\int_{\Omega_{\varepsilon}} \vert \varepsilon y+x_0\vert^l \vert \nabla z_{\varepsilon}(y) \vert^2dy} -  \lambda \varepsilon^{2} \displaystyle{\int_{\Omega_{\varepsilon}} w_{\varepsilon}(y) z_{\varepsilon}(y) dy}. $$

Then, as we can see the  three last terms have different weights and the exponents $k =2$ and $ l = 2$ are critical for our problem.

\bigskip
Let us finally point out that the lower dimension $ N= 3$ and the exponents $k<2$, $l< 2$  could also be interesting for this problem, but is not yet fully understood. For dimension $3$ in the scalar case we refer to \cite{BN}, \cite{CHL}.

The rest of the paper is organized in the following way : in Section~\ref{sec:main-results} we state our main result. In Section~3 we give a sufficient condition for the existence of minimizers. We give  precise estimates of the energy in Section~4. In Section~5 we discuss  the sign of minimizers. Section~6 contains some non-existence results.

\section{Main Results} \label{sec:main-results}
Let us suppose that there exists $ x_0 \in \Omega $ such that
$$ a(x_0) = \min \left\lbrace a(x), \ x \in \bar{\Omega} \right\rbrace.  $$
$$b(x_0) = \min \left\lbrace b(x), \ x \in \bar{\Omega} \right\rbrace.  $$
In this section, we will prove the existence of a solution to the problem \eqref{1.1}. For this we define:
\begin{equation} \label{S}
   S= \inf_{u \in H_0^1(\Omega), {\| u \|_{L^q}= 1}} \|\nabla u\|_2^{2} \, ,
\end{equation}
which corresponds to the best constant for the Sobolev embedding $H_{0}^{1}(\Omega)$ into $ L^{q}(\Omega). $

In this paper, we assume that $ a(x_0) = b(x_0)$, and we denote by $\gamma_0$ this value.

 We are now ready to state the main result of this paper
 \begin{theorem}  \label{theoreme1}
 Assume that $ a, b \in H^1(\Omega) \cap C(\bar{\Omega}) $ satisfy \eqref{1.3} and \eqref{1.4}, respectively. Let $ \lambda_{1}^{a}$, respectively  $ \lambda_{1}^{b}$,  be the first eigenvalue of $ -\div(a(x) \nabla \cdot)$, respectively $ -\div(b(x) \nabla \cdot)$ on $\Omega$ with zero Dirichlet boundary condition. Define $\tilde{\lambda}_1= \min
 \{\lambda_{1}^{a}, \lambda_{1}^{b}\}$. We have : \\
 $(1)$ If $N \geq 4 $, $k > 2$ and $l > 2, $ then $Q_\lambda $ is achieved for every $ \lambda \in \left( 0,  \lambda_{1}^{a} \right) \cap \left( 0,  \lambda_{1}^{b} \right)$. \\
  $(2)$ If $N \geq 5, $ $k = 2$ and $l = 2, $ then there exists a constant $ \gamma(N) = \frac{(N-2)N(N+2)}{8(N-1)}(A_2 + B_2) $ such that $ Q_\lambda$ is achieved for every $ \lambda \in \left( \gamma(N) ,  \tilde{\lambda}_{1} \right) $. \\
 $(3)$ If $N\geq 5, $ $k = 2$ and $l > 2, $ respectively $k > 2$ and $l = 2, $ then there exists a constant 
 $ \tilde{\gamma}(N) = \frac{N(N-2)(N+2)}{8(N-1)} $ such that  $Q_\lambda $ is achieved for every $ \lambda \in \left(\tilde{\gamma}(N)A_{2} ,  \tilde{\lambda}_{1} \right)$, respectively $ \lambda \in \left( \tilde{\gamma}(N)B_{2} ,  \tilde{\lambda}_{1} \right). $ \\
$(4)$ If $ N = 4, $ $k = 2$ and $l > 2$, respectively $k > 2$ and $l = 2, $ then there exists a minimizing solution of $ Q_\lambda$ for every $ \lambda \in \left( \tilde{A_2} , \tilde{\lambda}_{1} \right) $ respectively for every $ \lambda \in \left( \tilde{B_2} , \tilde{\lambda}_{1} \right)$ where $\tilde{A_2} $ and $\tilde{B_2} $ are some constants.
 \end{theorem}

We now proceed with the proof of Theorem \ref{theoreme1} as follows : First we show that \\ $ 0 \leq  Q_\lambda <\gamma_0 S $ and then we prove that this implies that the infimum $ Q_\lambda $ is achieved.

 \section{Sufficient Conditions for the Existence of Minimizers}

 We first prove the existence of $Q_\lambda, $ which is guaranteed by the following result :
 \begin{proposition} \mbox{} \label{WA}
 Let $ \varphi_{1}^{a}$, respectively $ \varphi_{1}^{b}$, be the first eigenfunction of $ -\div(a(x) \nabla \cdot)$, respectively $ -\div(b(x) \nabla \cdot) $, associated to the first eigenvalue $ \lambda_{1}^{a}$, respectively $ \lambda_{1}^{b}$.
We have
 \begin{itemize}
\item[(i)] Assume that  $0< \lambda < \tilde{\lambda}_1, $ then
 $$ Q_\lambda \geq 0. $$
 \item[(ii)] For $ \displaystyle{ \lambda \geq \frac{ \|\varphi_{1}^{a}  \|_{L^q} \| \varphi_{1}^{b} \|_{L^q}}{ \displaystyle{ \int_{\Omega} \varphi_{1}^{a} \varphi_{1}^{b} dx}} \vert \Omega \vert^{1-\frac{2}{q}} } \max( \lambda_{1}^{a}, \lambda_{1}^{b}), $ one has
$$ Q_\lambda \leq 0, $$
\end{itemize}
 \end{proposition}
 \textbf{Proof} $(i). $
 Let  $0< \lambda < \tilde{\lambda}_1(\Omega)$, and let $u  $ and $v$ be such that $\| u\|_{L^{q}} = \| v\|_{L^{q}} = 1. $ \\
 By the definitions of $  \lambda_{1}^{a}, \lambda_{1}^{b} $ and $ \tilde{\lambda}_1(\Omega)$  one has   \\
$$ E_\lambda(u, v) \geq \frac{\tilde{\lambda}_1}{2} \| u\|_{L^{2}}^{2} + \frac{\tilde{\lambda}_1}{2} \| v\|_{L^2}^{2} -  \lambda\int_{\Omega}  u v\  dx. $$
By applying the Cauchy-Schwarz inequality , we find
$$ E_\lambda(u, v) \geq \frac{\tilde{\lambda}_1}{2} \left[ \| u\|_{L^2}^{2} + \| v\|_{L^2}^{2} - 2 \| u\|_{L^2} \| v\|_{L^2} \right]. $$
Thus
$$ E_\lambda(u, v) \geq \frac{\tilde{\lambda}_1}{2} \left(  \| u\|_{L^2} -  \| v\|_{L^2} \right)^{2} \geq 0. $$
Consequently, $Q_\lambda \geq 0 $.

$(ii)$
We have
$$ Q_\lambda \leq E_\lambda(\frac{\varphi_{1}^{a}}{\|\varphi_{1}^{a}  \|_{L^q}}, \frac{\varphi_{1}^{b}}{\| \varphi_{1}^{b} \|_{L^q}}), $$
which implies that,
$$ Q_\lambda \leq \frac{1}{2 \|\varphi_{1}^{a}  \|_{L^q}^{2}} \int_{\Omega} a(x) \vert \nabla  \varphi_{1}^{a} \vert^2dx  + \displaystyle{ \frac{1}{2 \| \varphi_{1}^{b} \|_{L^q}^{2}} \int_{\Omega} b(x) \vert \nabla \varphi_{1}^{b} \vert^2dx } - \frac{\lambda}{ \|\varphi_{1}^{a}  \|_{L^q} \| \varphi_{1}^{b} \|_{L^q}} \displaystyle{\int_{\Omega} \varphi_{1}^{a} \varphi_{1}^{b} dx}. $$
 By the definitions of $  \varphi_{1}^{a}$ and $ \varphi_{1}^{b} ,$ one has
$$Q_\lambda \leq \frac{1}{2 \|\varphi_{1}^{a} \|_{L^q}^{2} } \lambda_{1}^{a} \int_{\Omega}  {\vert \varphi_{1}^{a} \vert}^2 dx + \frac{1}{2 \| \varphi_{1}^{b} \|_{L^q}^{2}} \lambda_{1}^{b} \int_{\Omega}  {\vert \varphi_{1}^{b} \vert}^2 dx  - \frac{\lambda}{ \|\varphi_{1}^{a}  \|_{L^q} \| \varphi_{1}^{b} \|_{L^q}} \int_{\Omega} \varphi_{1}^{a} \varphi_{1}^{b} dx. $$
Using the embedding of $ L^{q}$ into $ L^2$, there exists a positive constant $C_1 = \vert \Omega \vert^{\frac{1}{2}-{\frac{1}{q}}} $ such that
$$ Q_\lambda \leq \frac{1}{2 \|\varphi_{1}^{a} \|_{L^q}^{2}} \lambda_{1}^{a} C_1^{2} \|\varphi_{1}^{a} \|_{L^q}^{2} + \frac{1}{2 \| \varphi_{1}^{b} \|_{L^q}^{2} } \lambda_{1}^{b} C_1^{2} \| \varphi_{1}^{b} \|_{L^q}^{2} - \frac{\lambda}{ \|\varphi_{1}^{a}  \|_{L^q} \| \varphi_{1}^{b} \|_{L^q}} \int_{\Omega} \varphi_{1}^{a} \varphi_{1}^{b} dx. $$
Thus
$$ Q_\lambda \leq \vert \Omega \vert^{1-\frac{2}{q}} \max( \lambda_{1}^{a}, \lambda_{1}^{b}) -\frac{\lambda}{\|\varphi_{1}^{a}  \|_{L^q} \| \varphi_{1}^{b} \|_{L^q}}  \int_{\Omega} \varphi_{1}^{a} \varphi_{1}^{b} dx. $$
Therefore, $Q_\lambda\leq 0$, completing the proof of the proposition. \hfill $\blacksquare$
\begin{lemma} \label{F}
Let $0< \lambda < \tilde{\lambda}_1. $
If $ Q_\lambda < \gamma_0 S $, then the infimum in \eqref{1.1} is achieved.
  \end{lemma}
   \textbf{Proof} Let $\left\lbrace U_n \right\rbrace  \subset \left( H_{0}^{1}(\Omega)\right)^2 $ be a minimizing sequence for \eqref{1.1} that is,
 \begin{equation}  \label{ll}
 \| U_n \|_{L^{q}} = 1 \quad (\mbox{which means } \| u_n \|_{q} = 1, \quad \| v_n \|_{q} = 1).
 \end{equation}
\begin{equation}  \label{rr}
\displaystyle{ \frac{1}{2} \int_{\Omega} a(x) \vert \nabla u_n \vert^2dx + \frac{1}{2}\int_{\Omega} b(x) \vert \nabla v_n \vert^2dx -  \lambda \int_{\Omega}  u_n v_n  dx } =  Q_\lambda + o(1) \quad  \mbox{as} \quad  n \longrightarrow \infty.
\end{equation}
The sequence $\left\lbrace U_n \right\rbrace$ is bounded in $ \left( H_{0}^{1}(\Omega)\right)^2 $. Indeed, from \eqref{rr}, we have
 $$ \displaystyle{\frac{1}{2} \int_{\Omega} a(x) \vert \nabla u_n \vert^2dx + \frac{1}{2} \int_{\Omega} b(x) \vert \nabla v_n \vert^2 dx} = \displaystyle{  \lambda \int_{\Omega}  u_n  v_n  dx  + Q_\lambda + o(1). } $$
 Using H\"older's inequality we have
$$ \displaystyle{ \int_{\Omega} \vert u_n \vert \vert v_n \vert dx \leq  \| u_n \|_{L^{\frac{q}{q-1}}} \| v_n \|_{L^{q}}}, $$
then
$$ \displaystyle{\frac{1}{2} \left[  \int_{\Omega} a(x) \vert \nabla u_n \vert^2dx + \int_{\Omega} b(x) \vert \nabla v_n \vert^2 dx \right]
 \leq  \lambda  \| u_n \|_{L^{\frac{q}{q-1}}}  \| v_n \|_{L^{q}}  + Q_\lambda + o(1).} $$
 Consequently, we have
 $$ \| (u_n, v_n) \|_{\left( H_{0}^{1}(\Omega)\right)^2} \leq  C. $$
Then there is $ U= (u, v)$ , up to a subsequence, still denoted by $ U_n=(u_n, v_n)$ , such that

$$(u_n, v_n) \rightharpoonup (u, v) \quad \mbox{weakly in} \quad \left( H_{0}^{1}(\Omega)\right)^2, $$
$$(u_n, v_n) \longrightarrow (u, v) \quad \mbox{ strongly in} \quad  \left( L^2(\Omega)\right)^2, $$
$$(u_n, v_n) \longrightarrow (u, v)  \quad \mbox{a.e. on} \quad  \Omega. $$
with $ \| u \|_{L^{q}} \leq 1$ and $  \| v \|_{L^{q}} \leq 1. $ \\
Set $ w_n = u_n - u $ and $ z_n = v_n - v , $ so that
$$(w_n, z_n) \rightharpoonup (0, 0) \quad \mbox{weakly in} \quad \left( H_{0}^{1}(\Omega)\right)^2, $$
$$(w_n, z_n) \longrightarrow (0, 0) \quad \mbox{ strongly in} \quad  \left( L^2(\Omega)\right)^2, $$
$$(w_n, z_n) \longrightarrow (0, 0)  \quad \mbox{a.e. on} \quad  \Omega. $$
By \eqref{ll} and using the definition of $S$ and the fact that $\displaystyle{\min_{\bar{\Omega}} a(x)}= \gamma_0 > 0 $ and $\displaystyle{\min_{\bar{\Omega}} b(x)}= \gamma_0 > 0 , $ we have \\
$$  \displaystyle{\frac{1}{2} \left[  \int_{\Omega} a(x) \vert \nabla u_n \vert^2dx + \int_{\Omega} b(x) \vert \nabla v_n \vert^2 dx \right]  \geq \frac{1}{2} 2 \gamma_0 S.} $$
By \eqref{rr}, we have
$$ \displaystyle{\frac{1}{2} \left[  \int_{\Omega} a(x) \vert \nabla u_n \vert^2dx + \int_{\Omega} b(x) \vert \nabla v_n \vert^2dx \right] - \lambda \int_{\Omega}  u_n  v_n dx } =  Q_\lambda + o(1). $$
Which implies that
$$ \displaystyle{ \frac{1}{2} \left[ \int_{\Omega} a(x) \vert \nabla u_n \vert^2dx + \int_{\Omega} b(x) \vert \nabla v_n \vert^2 dx \right] } = \displaystyle{   Q_\lambda +  \lambda \int_{\Omega} u_n  v_n  dx  + o(1). } $$
Which gives
$$ \lambda \int_{\Omega}  u  v  dx \geq \gamma_0 S - Q_\lambda > 0 $$ and so $u \not\equiv 0 $ and $v \not\equiv 0. $ \\
Again using \eqref{rr} we obtain
$$ \displaystyle\frac{1}{2} \left[ \int_{\Omega} a(x) \vert \nabla u \vert^2dx + \int_{\Omega} a(x) \vert \nabla w_n \vert^2dx + \int_{\Omega} b(x) \vert \nabla v \vert^2dx\right.  + \left. \displaystyle \int_{\Omega} b(x) \vert \nabla z_n \vert^2 dx \right]$$
\begin{equation}  \label{aa}
- \lambda \int_{\Omega}  u_n  v_n dx  =  Q_\lambda + o(1).
\end{equation}

On the other hand, as $\left\lbrace w_n \right\rbrace $ and $\left\lbrace z_n \right\rbrace $  are bounded in $ L^{q}(\Omega)$ and $ w_n\longrightarrow 0$ for a.e. $x$ in $\Omega, $ $ z_n\longrightarrow 0$ for a.e. $x$ in $\Omega, $ then we can use the Brezis-Lieb lemma, see~\cite{BL}
$$ \| u+w_n \|_{L^{q}}^{q} = \| u \|_{L^{q}}^{q} + \| w_n \|_{L^{q}}^{q}+ o(1), $$
$$ \| v+z_n \|_{L^{q}}^{q} = \| v \|_{L^{q}}^{q} + \| z_n \|_{L^{q}}^{q}+ o(1). $$
Using \eqref{ll} we have
$$ 1 = \| u \|_{L^{q}}^{q} + \| w_n \|_{L^{q}}^{q} + o(1), $$
$$ 1 = \| v \|_{L^{q}}^{q} + \| z_n \|_{L^{q}}^{q} + o(1). $$
And so,
\begin{equation} \label{HA}
 1 \leq \| u \|_{L^{q}}^{2} + \| w_n \|_{L^{q}}^{2} + o(1),
 \end{equation}
\begin{equation}  \label{HB}
1 \leq \| v \|_{L^{q}}^{2} + \| z_n \|_{L^{q}}^{2} + o(1).
\end{equation}
We sum the equations \eqref{HA} and \eqref{HB} we obtain
\begin{equation} \label{2.7}
 2 \leq  \| u \|_{L^{q}}^{2} + \| v \|_{L^{q}}^{2} + \frac{1}{ \gamma_0 S} \left( \int_{\Omega} a(x) \vert \nabla w_n \vert^2dx + \int_{\Omega} b(x) \vert \nabla z_n \vert^2dx \right).
 \end{equation}
Since $ Q_\lambda > 0, $ this implies
\begin{equation} \label{ss}
 Q_\lambda \leq \frac{Q_\lambda}{2} \| u \|_{L^{q}}^{2} + \frac{Q_\lambda}{2} \| v \|_{L^{q}}^{2} + \frac{Q_\lambda}{ 2 \gamma_0 S} \left( \int_{\Omega} a(x) \vert \nabla w_n \vert^2dx + \int_{\Omega} b(x) \vert \nabla z_n \vert^2dx \right).
\end{equation}
Adding \eqref{aa} and \eqref{ss} we obtain
$$ \displaystyle{\frac{1}{2} \left[  \int_{\Omega} a(x) \vert \nabla u \vert^2dx + \int_{\Omega} a(x) \vert \nabla w_n \vert^2dx + \int_{\Omega} b(x) \vert \nabla v \vert^2dx + \int_{\Omega} b(x) \vert \nabla z_n \vert^2dx \right] } - \displaystyle{ \lambda \int_{\Omega}  u  v  dx } $$
$$ \leq  \frac{Q_\lambda}{2} \| u \|_{L^{q}}^{2} + \frac{Q_\lambda}{2} \| v \|_{L^{q}}^2 + \frac{Q_\lambda}{2 \gamma_0 S} \left( \int_{\Omega} a(x) \vert \nabla w_n \vert^2dx + \int_{\Omega} b(x) \vert \nabla z_n \vert^2dx \right). $$
Thus
$$ \displaystyle{ \frac{1}{2} \left[ \int_{\Omega} a(x) \vert \nabla u \vert^2dx + \int_{\Omega} b(x) \vert \nabla v \vert^2dx \right]  - \lambda \int_{\Omega}  u  v  dx } $$
$$\leq  \frac{Q_\lambda}{2} \| u \|_{L^{q}}^{2} +  \frac{Q_\lambda}{2} \| v \|_{L^{q}}^2 + \left[  \frac{Q_\lambda}{ 2 \gamma_0 S} - \frac{1}{2} \right] \left[ \int_{\Omega} a(x) \vert \nabla w_n \vert^2dx + \int_{\Omega} b(x) \vert \nabla z_n \vert^2 dx \right]. $$
Hence
$$ \displaystyle{ 2 E_\lambda(u,v) \leq Q_\lambda \| u \|_{L^{q}}^{2} + Q_\lambda \| v \|_{L^{q}}^{2} + \left[  \frac{Q_\lambda}{  \gamma_0 S} - 1 \right] \left[ \int_{\Omega} a(x) \vert \nabla w_n \vert^2dx + \int_{\Omega} b(x) \vert \nabla z_n \vert^2 dx \right]+ o(1)}, $$
and since $ Q_\lambda < \gamma_0 S$, we know that $u \neq 0$ and $v \neq 0. $
we deduce
$$ \displaystyle{ 2 E_\lambda\Big(\frac{u}{\| u \|_{L^q}},\frac{v}{\| v \|_{L^q}} \Big) \leq  \left( \frac{1}{\| u \|_{L^q}^{2}} - 1 \right)  \int_{\Omega} a(x) \vert \nabla u \vert^2 dx + \left( \frac{1}{\| v \|_{L^q}^{2}} - 1 \right)  \int_{\Omega} b(x) \vert \nabla v \vert^2 dx} $$
$$ \displaystyle{ +  Q_\lambda \| u \|_{L^{q}}^{2} + Q_\lambda \| u \|_{L^{q}}^{2} + \lambda \int_{\Omega}  u  v  dx - \lambda \int_{\Omega} \frac{ u  v}{\| u \|_{L^q} \| v \|_{L^q}}   dx. }$$
Then
$$ \displaystyle{ 2 E_\lambda(\frac{u}{\| u \|_{L^q}},\frac{v}{\| v \|_{L^q}} ) \leq  \left( \| u \|_{L^q}^{2} - 1 \right) \left[ Q_\lambda - \int_{\Omega} a(x)\frac{ \vert \nabla u \vert^2}{\| u \|_{L^q}^{2}} dx - 2 \lambda \int_{\Omega} \frac{ u  v}{\| u \|_{L^q} \| v \|_{L^q}} dx \right] }$$
$$ \displaystyle{ + \left( \| v \|_{L^q}^{2} - 1 \right)   \left[ Q_\lambda - \int_{\Omega} b(x)\frac{ \vert \nabla v \vert^2}{\| v \|_{L^q}^{2}} dx - 2 \lambda \int_{\Omega} \frac{ u  v}{\| u \|_{L^q} \| v \|_{L^q}} dx \right]}$$
$$ \displaystyle{ + \lambda \int_{\Omega}  u  v  dx - \lambda \int_{\Omega} \frac{ u  v}{\| u \|_{L^q} \| v \|_{L^q}}   dx + 2 \lambda \left( \| u \|_{q}^{2} - 1 \right) \int_{\Omega} \frac{ u  v}{\| u \|_{L^q} \| v \|_{L^q}} dx }$$
$$ + 2 \lambda \left( \| v \|_{L^q}^{2} - 1 \right) \int_{\Omega} \frac{ u  v}{\| u \|_{L^q} \| v \|_{L^q}} dx + 2 Q_\lambda. $$
On one hand, we have $ \| u \|_{L^q}^{2} - 1 \leq 0 $ and $ \| v \|_{L^q}^{2} - 1 \leq 0, $ we obtain
$$ \displaystyle{ 2 E_\lambda\Big(\frac{u}{\| u \|_{L^q}},\frac{v}{\| v \|_{L^q}} \Big) \leq \lambda \left[ \frac{\left( \| u \|_{L^q} \| v \|_{L^q}- 1\right) + 2 \left( \| u \|_{L^q}^{2} - 1 \right) + 2 \left( \| v \|_{L^q}^{2} - 1 \right)}{\| u \|_{L^q} \| v \|_{L^q}} \right] \int_{\Omega}  u  v dx + 2 Q_\lambda. }$$
On the other hand, we have $\| u \|_{L^q} \| v \|_{L^q}- 1 \leq 0, $ then
 $$ \displaystyle{ 2 E_\lambda\left(\frac{u}{\| u \|_{L^q}},\frac{v}{\| v \|_{L^q}} \right) \leq 2 Q_\lambda. } $$

This means that $(u,v)$ is a minimum of $Q_\lambda.$
\hfill $\blacksquare$
 \section{Precise estimates of the energy}

 \begin{proposition}  \mbox{} \label{Q}
 \begin{itemize}
\item[(a)] For $N \geq 4,$ and $ k> 2, $ $ l > 2$  we have
$$ Q_\lambda < \gamma_0 S, \quad \mbox{for all} \quad \lambda > 0. $$
\item[(b)] For $N = 4,$ and $ k = 2, $ $ l = 2$  we have
$$ Q_\lambda <  \gamma_0 S, \quad \mbox{for all} \quad \lambda > A_2 + B_2. $$
\item[(c)] For $N \geq 5, $ and $ k = 2, $ $ l = 2$  we have
$$ Q_\lambda <  \gamma_0 S, \quad \mbox{for all} \quad \lambda > m_N \left( A_2+ B_2 \right).  $$
\item[(d)]  For $N \geq 5,$ and $ k = 2, $ $ l > 2$ (respectively $ k > 2, $ $ l = 2$)  we have
 $$ Q_\lambda <  \gamma_0 S, \quad \mbox{for all} \quad \lambda > m_N A_2 \quad (\mbox{respectively} \quad \lambda > m_N B_2 ). $$
\item[(e)] For $N = 4,$ and $ k = 2, $ $ l > 2$, respectively  $ k > 2, $ $ l = 2$ we have
 $$ Q_\lambda <  \gamma_0 S, \quad \mbox{for all} \quad \lambda >  A_2,\,\,\, \textrm{respectively\,} \mbox{for all} \quad \lambda >  B_2. $$
 \end{itemize}
 where $A_2, B_2$ are defined by \eqref{1.5} and \eqref{1.6}, $ K_3 = \displaystyle{\int_{\mathbb{R}^N} \frac{1}{\left( 1+{\vert x \vert}^2 \right)^{N-2}}dx } $ and \\ $m_N=\frac{N(N-2)(N+2)}{8(N-1)}. $
\end{proposition}
\noindent \textbf{Proof}  We shall estimate the ratio $E_\lambda(u, v)$ defined in \eqref{1.1}, with $u=u_{x_0,\varepsilon}=\zeta  U_{\varepsilon}(x-x_0), $  where for $ x \in \mathbb{R}^N, \quad \displaystyle{  U_{\varepsilon}(x)= \frac{\varepsilon^{\frac{N-2}{4}}}{\left( \varepsilon + {\vert x \vert}^2 \right)^{\frac{N-2}{2}}} } $ and $ \zeta \in C_c^\infty{(\Omega)} $ with $ \zeta \geq 0 $ and $ \varphi \equiv 1$ on a neighborhood of $x_0; $ for more details see  \cite{BN} \cite{T}.  \\
  We recall from \cite{BN} that
\begin{equation} \label{2.9}
  \int_{\Omega} {\vert \nabla u_{x_0, \varepsilon} (x) \vert}^2 dx  = K_1  + O(\varepsilon^{\frac{N-2}{2}}),
 \end{equation}
\begin{equation} \label{2.10}
 \left(  \int_{\Omega} {\vert u_{x_0, \varepsilon} (x) \vert}^{q} dx \right)^{\frac{2}{q}}  = K_2 + O(\varepsilon^{\frac{N-2}{2}}),
\end{equation}
 \begin{equation} \label{2.11}
 \displaystyle{ \int_{\Omega} \vert u_{x_0, \varepsilon}(x) \vert^2  dx = \left\{
   \begin{array}{ccl}
    K_3 \varepsilon +  O(\varepsilon^{\frac{N-2}{2}}),  & \mbox{ if } & N \geqslant 5,\\
   \frac{\omega_4}{2} \varepsilon \vert \ln \varepsilon \vert +  o(\varepsilon \vert \ln \varepsilon \vert), &\mbox{ if }&  N=4, \\
      \end{array}
    \right.}
 \end{equation}
  where $K_1$ and $K_2$ are positive constants with $\frac{K_1}{K_2} = S $, $w_4$ is a the area of $S^3$  \\
    and $K_3 = \displaystyle{\int_{\mathbb{R}^N} \frac{1}{\left( 1+{\vert x \vert}^2 \right)^{N-2}}dx }. $ \\
    We know by \cite{HY}, the following estimates :
   $$ \frac{\varepsilon^{\frac{N-2}{2}}}{2} \int_{\Omega} a(x) {\vert \nabla u_{x_0, \varepsilon} (x) \vert}^2 dx $$
 \begin{equation}  \label{2.13}
 \leq\left\{
       \begin{array}{lclll}
    \frac{a(x_0) K_1}{2} +O(\varepsilon^{\frac{N-2}{2}})\quad& \mbox{if}\ N\geq4 \ \mbox{and}\ N-2<k , \\
  \frac{a(x_0) K_1}{2} + \frac{C_k}{2} \varepsilon^{\frac{k}{2}}+ o(\varepsilon^{\frac{k}{2}}) \quad & \mbox{if}\ N\geq4 \ \mbox{and}\ N-2>k , \\
  \frac{a(x_0) K_1}{2}+ \frac{(N-2)^2 \omega_N(A_{N-2} + M)}{2} \varepsilon^{\frac{N-2}{4}} \vert \ln \varepsilon \vert +o(\varepsilon^{\frac{N-2}{2}} \vert \ln \varepsilon \vert) \quad & \mbox{if}\ N>4\ \mbox{and}\ k=N-2 ,    \\
 \frac{a(x_0) K_1}{2} + A_2 \omega_4 \varepsilon \vert \ln \varepsilon \vert + o(\varepsilon \vert \ln \varepsilon \vert) \quad &\mbox{if}\ N=4\ \mbox{and}\ k=2 \\
      \end{array}
    \right.
    \end{equation}
    and
    $$  \varepsilon^{\frac{N-2}{2}} \int_{\Omega} b(x) {\vert \nabla u_{x_0, \varepsilon} (x) \vert}^2 dx $$
    \begin{equation} \label{2.14}
 \leq\left\{
       \begin{array}{lclll}
    \frac{b(x_0) K_1}{2}+O(\varepsilon^{\frac{N-2}{2}})\quad& \mbox{if}\ N\geq4 \ \mbox{and}\ N-2<l , \\
 \frac{ b(x_0) K_1}{2} + \frac{D_l}{2} \varepsilon^{\frac{l}{2}}+ o(\varepsilon^{\frac{l}{2}}) \quad & \mbox{if}\ N\geq4 \ \mbox{and}\ N-2>l , \\
 \frac{b(x_0) K_1}{2}+ \frac{(N-2)^2 \omega_N(B_{N-2} + M)}{4} \varepsilon^{\frac{N-2}{2}} \vert \ln \varepsilon \vert +o(\varepsilon^{\frac{N-2}{2}} \vert \ln \varepsilon \vert) \quad & \mbox{if}\ N>4\ \mbox{and}\ l=N-2 ,    \\
 \frac{b(x_0) K_1}{2} + B_2 \omega_4 \varepsilon \vert \ln \varepsilon \vert + o(\varepsilon \vert \ln \varepsilon \vert) \quad &\mbox{if}\ N=4 \ \mbox{and}\ l = 2 \\
      \end{array}
    \right.
    \end{equation}
    Where $ K_1 = \displaystyle{ (N-2)^2 \int_{\mathbb{R}^N} \frac{{\vert y \vert}^2}{\left( 1+{\vert y \vert}^2 \right)^{N}}dy }, $ $K_2 = \displaystyle{ \left( \int_{\mathbb{R}^N} \frac{dy}{\left( 1+{\vert y \vert}^2 \right)^{N}}\right)^{\frac{N-2}{N}}  }, $\\
     $ C_k = \displaystyle{ (N-2)^2 A_k \int_{\mathbb{R}^N} \frac{{\vert y \vert}^{k+2}}{\left( 1+{\vert y \vert}^2 \right)^{N}}dy }, $ $  D_k = \displaystyle{ (N-2)^2 B_k \int_{\mathbb{R}^N} \frac{{\vert y \vert}^{l+2}}{\left( 1+{\vert y \vert}^2 \right)^{N}}dy }, $ $ M $ and $ M' $ are positive constants. \\
 We claim that, as $\varepsilon\rightarrow0, $  by adding \eqref{2.10}, \eqref{2.11}, \eqref{2.13} and \eqref{2.14}
 we obtain
  \noindent $E(u_{x_0, \varepsilon}) \leq$
  $$\displaystyle{ \hspace{-2cm}\qquad
\left\{
\begin{array}{lclll}
\frac{a(x_0) K_1+o(\varepsilon)}{2 K_2} + \frac{b(x_0) K_1 + o(\varepsilon)}{ 2 K_2} - \lambda \frac{K_3}{K_2} \varepsilon  \quad&\mbox{if} \left\{
\begin{array}{lcl}
  N \geq 5 \\
k> 2, l> 2,
\end{array}
\right. \\
\frac{a(x_0) K_1+ C_2 \varepsilon + o(\varepsilon)}{ 2 K_2} + \frac{b(x_0) K_1+ D_2 \varepsilon + o(\varepsilon)}{ 2 K_2} - \lambda \frac{K_3}{K_2} \varepsilon  \quad&\mbox{if} \left\{
\begin{array}{lcl}
  N \geq 5 \\
k = 2, l = 2,
\end{array}
\right. \\
\frac{a(x_0) K_1 + o(\varepsilon)}{ 2 K_2} + \frac{b(x_0) K_1+o(\varepsilon)}{ 2 K_2} -  \lambda \varepsilon \frac{\omega_4}{2} \frac{\vert \log \varepsilon \vert}{K_2} + o( \varepsilon \vert \log \varepsilon \vert)  \quad&\mbox{if} \left\{
\begin{array}{lcl}
  N = 4 \\
k> 2, l> 2,
\end{array}
\right. \\
\frac{a(x_0) K_1 + 2 \omega_4 A_2 \varepsilon \vert \log \varepsilon \vert + o(\vert \log \varepsilon \vert)}{ 2 K_2} + \frac{b(x_0) K_1+ 2 \omega_4 B_2 \varepsilon \vert \log \varepsilon \vert  + o(\varepsilon \vert \log \varepsilon \vert)}{ 2 K_2} -  \lambda \varepsilon \frac{\omega_4}{2 K_2} \vert \log \varepsilon \vert  \quad&\mbox{if} \left\{
\begin{array}{lcl}
  N = 4 \\
k = 2, l = 2,
\end{array}
\right. \\
\frac{a(x_0) K_1+ C_2 \varepsilon + o(\varepsilon)}{ 2 K_2} + \frac{b(x_0) K_1 + o(\varepsilon)}{K_2} -  \lambda \frac{K_3}{ K_2} \varepsilon + O(1) \quad&\mbox{if} \left\{
\begin{array}{lcl}
  N \geq 5 \\
k = 2, l > 2,
\end{array}
\right. \\
\frac{a(x_0) K_1 + o(\varepsilon)}{2 K_2} + \frac{b(x_0) K_1+ D_2 \varepsilon + o(\varepsilon)}{ 2 K_2} - \lambda \frac{K_3}{K_2} \varepsilon + O(1)  \quad&\mbox{if} \left\{
\begin{array}{lcl}
  N \geq 5 \\
k > 2, l = 2,
\end{array}
\right. \\
\frac{a(x_0) K_1 + o( \varepsilon )}{2 K_2} + \frac{b(x_0) K_1+ 2 \omega_4 B_2 \varepsilon \vert \log \varepsilon \vert  + o(\varepsilon \vert \log \varepsilon \vert)}{2 K_2} - \lambda \varepsilon \frac{\omega_4}{2 K_2} \vert \log \varepsilon \vert + o(\varepsilon \vert \log \varepsilon \vert)  \quad&\mbox{if} \left\{
\begin{array}{lcl}
  N = 4 \\
k > 2, l = 2,
\end{array}
\right. \\
\frac{a(x_0) K_1+ 2 \omega_4 A_2 \varepsilon \vert \log \varepsilon \vert + o(\vert \log \varepsilon \vert)}{2 K_2} + \frac{b(x_0) K_1 + o(\varepsilon)}{2 K_2} - \lambda \varepsilon \frac{\omega_4}{2 K_2} \vert \log \varepsilon \vert + o(\varepsilon \vert \log \varepsilon \vert) \quad&\mbox{if} \left\{
\begin{array}{lcl}
  N = 4 \\
k = 2, l > 2, \end{array}
\right. \\
 \frac{a(x_0) K_1+ C_k {\varepsilon}^{\frac{k}{2}} + o({\varepsilon}^{\frac{k}{2}})}{2 K_2} + \frac{b(x_0) K_1 + 2 \omega_4 B_2 \varepsilon \vert \log \varepsilon \vert  + o(\varepsilon \vert \log \varepsilon \vert)}{ 2 K_2} -  \lambda \varepsilon \frac{\omega_4}{2 K_2} \vert \log \varepsilon \vert + o(\varepsilon \vert \log \varepsilon \vert)   \quad&\mbox{if} \left\{
\begin{array}{lcl}
  N = 4 \\
k < 2, l = 2,
\end{array}
\right. \\
\frac{a(x_0) K_1+ 2 \omega_4 A_2 \varepsilon \vert \log \varepsilon \vert + o( \varepsilon \vert \log \varepsilon \vert)}{ 2 K_2} + \frac{b(x_0) K_1+  D_k {\varepsilon}^{\frac{l}{2}} + o({\varepsilon}^{\frac{l}{2}})}{ 2 K_2} - \lambda \varepsilon \frac{\omega_4}{2 } \vert \log \varepsilon \vert + o(\varepsilon \vert \log \varepsilon \vert) \quad&\mbox{if} \left\{
\begin{array}{lcl}
  N = 4 \\
k = 2, l < 2,
\end{array}
\right. \\
\end{array}
\right.
}$$

\noindent Thus

\noindent $ Q_{\lambda}  \leq   $
$$
\qquad \left\{
\begin{array}{lclll}
(\frac{a(x_0) + b(x_0)}{2})S - \lambda \frac{K_3}{K_2} \varepsilon + o(\varepsilon)  \quad&\mbox{if} \left\{
\begin{array}{lcl}
  N \geq 5 \\
k> 2, l> 2,
\end{array}
\right. \\
(\frac{a(x_0) + b(x_0)}{2})S - (\lambda - \frac{C_2}{2 K_3} - \frac{D_2}{2 K_3} ) \frac{K_3}{K_2} \varepsilon + o(\varepsilon)\quad&\mbox{if} \left\{
\begin{array}{lcl}
  N \geq 5 \\
k = 2, l = 2,
\end{array}
\right. \\
(\frac{a(x_0) + b(x_0)}{2})S - \lambda \frac{\omega_4}{K_2} \varepsilon \vert \log \varepsilon \vert + o( \varepsilon \vert \log \varepsilon \vert) \quad&\mbox{if} \left\{
\begin{array}{lcl}
  N = 4 \\
k > 2, l > 2,
\end{array}
\right. \\
(\frac{a(x_0) + b(x_0)}{2})S - \frac{\omega_4}{K_2}( \lambda - (A_2 + B_2)) \varepsilon \vert \log \varepsilon \vert + o( \varepsilon \vert \log \varepsilon \vert) \quad& \mbox{if} \left\{
\begin{array}{lcl}
  N = 4 \\
k = 2, l = 2,
\end{array}
\right. \\
 (\frac{a(x_0) + b(x_0)}{2})S - ( \lambda - \frac{C_2}{2 K_3}) \frac{K_3}{K_2} \varepsilon + o(\varepsilon)\quad&\mbox{if} \left\{
\begin{array}{lcl}
  N \geq 5 \\
k = 2, l > 2,
\end{array}
\right. \\
(\frac{a(x_0) + b(x_0)}{2})S - ( \lambda - \frac{D_2}{2 K_3}) \frac{K_3}{K_2} \varepsilon + o(\varepsilon)\quad&\mbox{if} \left\{
\begin{array}{lcl}
  N \geq 5 \\
k > 2, l = 2,
\end{array}
\right. \\
(\frac{a(x_0) + b(x_0)}{2})S - \frac{\omega_4}{K_2}\left[  \lambda - B_2 \right]  \varepsilon \vert \log \varepsilon \vert + o( \varepsilon \vert \log \varepsilon \vert) \quad&\mbox{if} \left\{
\begin{array}{lcl}
  N = 4 \\
k > 2, l = 2,
\end{array}
\right. \\
(\frac{a(x_0) + b(x_0)}{2})S - \frac{\omega_4}{K_2}\left[  \lambda - A_2 \right]  \varepsilon \vert \log \varepsilon \vert + o( \varepsilon \vert \log \varepsilon \vert) \quad&\mbox{if} \left\{
\begin{array}{lcl}
  N = 4 \\
k = 2, l > 2,
\end{array}
\right. \\
(\frac{a(x_0) + b(x_0)}{2})S - \frac{\omega_4}{K_2}\left[  \lambda - B_2 \right] \varepsilon \vert \log \varepsilon \vert + \frac{C_k}{K_2} {\varepsilon}^{\frac{k}{2}} + o( \varepsilon \vert \log \varepsilon \vert) +o({\varepsilon}^{\frac{k}{2}}) \quad&\mbox{if} \left\{
\begin{array}{lcl}
  N = 4 \\
k < 2, l = 2,
\end{array}
\right. \\
(\frac{a(x_0) + b(x_0)}{2})S - \frac{\omega_4}{K_2}\left[  \lambda - A_2 \right] \varepsilon \vert \log \varepsilon \vert + \frac{D_k}{K_2} {\varepsilon}^{\frac{l}{2}} +o({\varepsilon}^{\frac{l}{2}})+ o( \varepsilon \vert \log \varepsilon \vert)  \quad&\mbox{if} \left\{
\begin{array}{lcl}
  N = 4 \\
k = 2, l < 2,
\end{array}
\right. \\
 \end{array}
\right.
$$
 where $ \frac{C_2}{K_3} = \frac{N(N-2)(N+2)}{4(N-1)} A_2 $ \quad and \quad $ \frac{D_2}{K_3} = \frac{N(N-2)(N+2)}{4(N-1)} B_2.$  \\
From these estimates we get the desired result. \hfill $\blacksquare$

Combining Lemma \ref{F} and Proposition \ref{Q} we conclude that the infimum in \eqref{1.1} is achieved. This leads to the conclusion of Theorem \ref{theoreme1}.

\begin{remark}
We  do not know if  we can go strictly  below $\frac{a(x_0) + b(x_0)}{2}S$ in the case where $ N \geq 4$, $k < 2$ and $ l < 2$. 
\end{remark}
\section{The Sign of the Minimizers}
Now, we will discuss the positivity of the solutions.
\begin{proposition} \mbox{}
 \begin{itemize}
 \item[(i)] If $ Q_\lambda $
 is achieved in $(u, v) $ then $uv\geq 0. $ 
  \item[(ii)] There exists a solutions $u$ and  $v$ of the minimization problem \eqref{1.1} such that $ u \geq 0$ and $v \geq 0. $
 \end{itemize}
 \end{proposition}  \label{pr1}

 \textbf{Proof} (i)
Set
$ \displaystyle{F(u, v)= \frac{1}{2} \int_{\Omega} a(x) \vert \nabla u \vert^2dx + \frac{1}{2} \int_{\Omega} b(x) \vert \nabla v \vert^2dx -  \lambda \displaystyle{\int_{\Omega} u(x)v(x) dx}}. $  \\ The inequality
  $$ F(u, v) \leq F(\vert u \vert, \vert v \vert) $$ gives that $ - \displaystyle{\lambda \int_{\Omega} u v dx  \leq - \lambda \int_{\Omega} \vert u v } \vert dx  $,  then
\begin{equation} \label{2.15}
  \displaystyle{ \int_{\Omega} ( u(x) v(x) - \vert u(x) v(x) \vert )dx \geq 0. }
\end{equation}
 We have always  $ u(x)v(x) - \vert u(x) v(x) \vert \leq 0 $ $ \forall x \in \Omega $. Hence
 \begin{equation}  \label{2.16}
  \int_{\Omega} ( u(x) v(x) - \vert u(x) v(x) \vert )dx \leq 0.
  \end{equation}
Combining \eqref{2.15} and \eqref{2.16} we obtain that $ uv - \vert u v \vert = 0 $ and hence $ \vert u v \vert = u v. $ \\
 Finally
for all $x \in \Omega, u(x) v(x) \geq 0. $ \\

 (ii) From Lemma \ref{F} we know there exists a minimum  $(u, v).$ By Proposition~\ref{pr1}~(i) we know that $\vert u v \vert = u v. $ \\
 On the other hand, we have
 $$ F(\vert u_n \vert, \vert v_n \vert)= \frac{1}{2} \int_{\Omega} a(x) \vert \nabla u_n \vert^2dx + \frac{1}{2} \int_{\Omega} b(x) \vert \nabla v_n \vert^2dx -  \lambda \displaystyle{\int_{\Omega} \vert u_n v_n \vert dx}  = F(u_n, v_n). $$
 Then if we have $(u_n, v_n)$ are solutions of \eqref{1.1} then $( \vert u_n \vert, \vert v_n \vert)$ are also minimizing solutions.

 Thus, when dealing with \eqref{1.1}, one can assume without loss of generality that $u \geq 0$ and $v \geq 0.$  \hfill $\blacksquare$

 \begin{remark}
 Since $u$ and $v$ are  minimizers for \eqref{1.1}  we obtain Lagrange multipliers $\Lambda_1, \Lambda_2 \in \mathbb{R}$ such that the Euler-Lagrange equation formula associated to \eqref{1.1} is
\begin{equation} \label{BN}
 \left\{
       \begin{array}{rclll}
 - \div(a(x) \nabla u) - \lambda v  &=& \Lambda_1 u^{2^{*}-1}  &\mbox{in} & \Omega \\
 - \div(b(x) \nabla v) - \lambda u  &= & \Lambda_2  v^{2^{*}-1} &\mbox{in} &\Omega \\
 \vert\vert u\vert\vert_{2^*} = \vert\vert u\vert\vert_{2^*} &=& 1\\
  u \geq 0, \quad   v &\geq& 0 &\mbox{ in } &\Omega \\
  u =v&=& 0  &\mbox{on}  &\partial \Omega,\\
 \end{array}
    \right.
\end{equation}
 such that $ \frac{\Lambda_1 + \Lambda_1}{2} = Q_\lambda \geq 0 $  according to Proposition \ref{WA}. \\
 \end{remark}

\section{Non-existence Results}
In this section, we assume that $a, b$ satisfy \eqref{1.3} and \eqref{1.4} with $ k>2$, $ l>2 $ and \eqref{1.5},  \eqref{1.6} if  $ k\leq 2$, $ l \leq 2$. we obtain in these cases  a few non-existence results.  We define
 $$  \omega(a, b):= \inf_{(u, v) \in \left( H_0^1(\Omega) \right)^2 \setminus \{0 \}} \phi_{a, b} (u, v) $$
where
$$ \phi_{a, b} (u, v):= \frac{1}{4} \frac{\displaystyle{\int_{\Omega} \left( \tilde{a}(x) {\vert \nabla u \vert}^2 + \tilde{b}(x) {\vert \nabla v \vert}^2\right) dx}}{ \displaystyle{\int_{\Omega} u(x) v(x) dx }}. $$
We see that $ \omega(a, b) \in \left[ -\infty, + \infty \right[. $

The main goal of this section is the non-existence results.
\begin{theorem} \label{XX}
 Assume that $\lambda \leq \omega(a, b)$ and $\Omega$ is a strictly star-shaped domain with respect to $x_0$. Then \eqref{1.1} has no minimizers.
 \end{theorem}
 The proof of this result  follows from the following Pohozaev identity; see \cite{P}.

  \begin{proposition}\label{YY}
 If $(u, v)$ is a solution of \eqref{BN} then $(u, v)$ satisfies the following identity :
 $$ \quad 2 \lambda \int_{\Omega} u(x) v(x) dx - \frac{1}{2} \int_{\Omega} \nabla b(x) \cdot (x-x_0) {\vert \nabla v(x) \vert}^2 dx -\frac{1}{2} \int_{\Omega} \left[  \nabla a(x) \cdot (x-x_0)\right]{\vert \nabla u(x)  \vert}^2 dx $$
$$ = \frac{1}{2} \int_{\partial \Omega} a(x) \left[ (x-x_0)\cdot \mathbf{n} \right]  {\Big\vert \frac{\partial u}{\partial \nu} \Big\vert}^2 d\sigma(x) + \frac{1}{2} \int_{\partial \Omega} b(x) \left[ (x-x_0)\cdot \mathbf{n} \right] {\Big\vert \frac{\partial v}{\partial \nu} \Big \vert}^2 d\sigma(x). $$
where $\mathbf{n}$ denotes the outward normal to $\partial \Omega. $
 \end{proposition}
 \textbf{Proof}
  Suppose that $(u, v)$ is a solution of Problem~\eqref{BN}. We multiply the first equation in the system by $ \nabla u(x) \cdot (x-x_0) $ and we integrate by parts, leading to 
 \begin{equation} \label{3.2}
  \hspace{-7cm} \Lambda_1 \int_{\Omega} \vert u (x)\vert^{q-2} u(x) \nabla u(x) \cdot (x-x_0) dx = - \frac{N}{q}  \Lambda_1,
   \end{equation}
   
  \begin{equation}   \label{3.3}
   \lambda  \int_{\Omega} v(x) \nabla u(x) \cdot (x-x_0) dx = - \lambda  \int_{\Omega} u(x) \nabla v(x) \cdot (x-x_0) dx - N \lambda \int_{\Omega} v(x) u(x) dx
   \end{equation}
   and
  \begin{align} 
 \int_{\Omega} - \div(a(x) \nabla u(x)) \nabla u(x)\cdot(x-x_0) dx =& - \frac{N-2}{2} \int_{\Omega} a(x) {\vert \nabla u(x) \vert}^2 dx \nonumber \\
& - \frac{1}{2} \int_{\Omega} \nabla a(x)\cdot(x-x_0) {\vert \nabla u(x)\vert}^2 dx \nonumber  \\
& - \frac{1}{2} \int_{\partial \Omega} a(x) {\Big\vert \frac{\partial u}{\partial \nu} \Big\vert}^2 \left( (x-x_0)\cdot \mathbf{n} \right)  d\sigma(x), \label{3.1}
\end{align}
 where $\mathbf{n}$ denotes the outward normal to $\partial \Omega. $

Combining \eqref{3.2}, \eqref{3.3}and  \eqref{3.1} we get
  $$
 - \frac{N-2}{2} \int_{\Omega} a(x) {\vert \nabla u(x)\vert}^2 dx - \frac{1}{2} \int_{\Omega} \nabla a(x) \cdot (x-x_0) {\vert \nabla u(x) \vert}^2 dx - \frac{1}{2} \int_{\partial \Omega} a(x) {\Big\vert \frac{\partial u}{\partial \nu} \Big\vert}^2 \left( (x-x_0) \cdot \mathbf{n}\right) d\sigma(x) $$
  \begin{equation} \label{3.4}
 =  \lambda  \int_{\Omega} v(x) \nabla u(x) \cdot (x-x_0) dx - \frac{N}{q}  \Lambda_1.
 \end{equation}
 Similarly, we multiply the second equation of \eqref{BN} by $ \nabla v(x) \cdot (x-x_0) $ and we integrate by parts, we get
 $$ -\frac{N-2}{2} \int_{\Omega} b(x) {\vert \nabla v(x) \vert}^2 dx - \frac{1}{2} \int_{\Omega} \nabla b(x) \cdot (x-x_0) {\vert \nabla v(x) \vert}^2 dx - \frac{1}{2} \int_{\partial \Omega} b(x) {\Big\vert \frac{\partial v}{\partial \nu} \Big\vert}^2 \left( (x-x_0)\cdot \mathbf{n} \right) d\sigma(x) $$
  \begin{equation} \label{3.5}
 =   \lambda  \int_{\Omega} u(x) \nabla v(x) \cdot (x-x_0) dx  - \frac{N}{q}  \Lambda_2.
 \end{equation}
 On the other hand, multiplying the equations in \eqref{BN} by $\frac{N-2}{2} u$ and  $\frac{N-2}{2} v $ respectively, and integrating and summing the obtained results, we get

 \begin{equation}  \label{3.6}
 \frac{N-2}{2} \int_{\Omega} a(x) {\vert \nabla u(x) \vert}^2 dx = \frac{N-2}{2} \lambda\int_{\Omega} v(x) u(x) dx  + \Lambda_1 \frac{N-2}{2}
 \end{equation}
 and
 \begin{equation}  \label{3.7}
 \frac{N-2}{2} \int_{\Omega} b(x) {\vert \nabla v(x) \vert}^2 dx = \frac{N-2}{2} \lambda\int_{\Omega} u(x) v(x)  dx + \Lambda_2 \frac{N-2}{2}.
 \end{equation}
Combining \eqref{3.4}, and \eqref{3.6} we obtain

\begin{align} \label{3.8}
 & - \left( \frac{N- 2}{2} \right)  \lambda \int_{\Omega} u(x) v(x) dx - \lambda  \int_{\Omega} v(x) \left[ \nabla u(x) \cdot (x-x_0) \right] dx \nonumber \\
  &= \frac{1}{2} \int_{\Omega} {\vert \nabla u \vert}^2 \left[  \nabla a(x) \cdot (x-x_0)\right] dx  + \frac{1}{2} \int_{\partial \Omega} a(x) {\vert \frac{\partial u}{\partial \nu} \vert}^2 \left[ (x-x_0)\cdot \mathbf{n} \right] d\sigma(x).
 \end{align}
 On the other hand, combining \eqref{3.5}, and \eqref{3.7} we get
\begin{align} \label{3.9}
 & - \left( \frac{N -2}{2} \right)  \lambda \int_{\Omega} u(x) v(x) dx - \lambda  \int_{\Omega} u(x) \left[ \nabla v(x) \cdot (x-x_0) \right] dx \nonumber \\
  &= \frac{1}{2} \int_{\Omega} {\vert \nabla v(x) \vert}^2 \left[  \nabla b(x) \cdot (x-x_0)\right] dx  + \frac{1}{2}
   \int_{\partial \Omega} b(x) {\vert \frac{\partial v}{\partial \nu} \vert}^2 \left[ (x-x_0)\cdot \mathbf{n} \right] d\sigma(x).
 \end{align}

Adding  \eqref{3.9} in \eqref{3.8}  and using \eqref{3.3} we find
$$ \quad 2 \lambda \int_{\Omega} u(x) v(x) dx - \frac{1}{2} \int_{\Omega} \nabla b(x) \cdot (x-x_0) {\vert \nabla v(x) \vert}^2 dx -\frac{1}{2} \int_{\Omega} \left[  \nabla a(x) \cdot (x-x_0)\right]{\vert \nabla u \vert}^2 dx $$
$$ = \frac{1}{2} \int_{\partial \Omega} a(x) {\vert \frac{\partial u}{\partial \nu} \vert}^2  \left[ (x-x_0)\cdot \mathbf{n} \right]  d\sigma(x) + \frac{1}{2} \int_{\partial \Omega} b(x) {\vert \frac{\partial v}{\partial \nu} \vert}^2 \left[ (x-x_0)\cdot \mathbf{n} \right]  d\sigma(x). $$
Finally we get the Proposition. \hfill $\blacksquare$

Now let us prove Theorem \ref{XX}.  Assume that $\lambda \leq \omega(a, b)$ and let $\Omega$ be strictly starshaped domain with respect to $x_0$ then
$$ (x-x_0)\cdot \mathbf{n} > 0, \qquad \mbox{for all} \  x \in \partial \Omega. $$
Suppose that $(u, v)$ is a solution of \eqref{BN}.
By the Proposition \ref{YY} we get
$$ 2 \lambda \int_{\Omega} u(x) v(x) dx - \frac{1}{2} \int_{\Omega} \nabla b(x) \cdot (x-x_0) {\vert \nabla v(x) \vert}^2 dx -\frac{1}{2} \int_{\Omega} \left[  \nabla a(x) \cdot (x-x_0)\right]{\vert \nabla u(x) \vert}^2 dx > 0. $$
It follows that
$$ \displaystyle{ \lambda > \frac{1}{4} } \displaystyle{ \inf_{ (u, v) \in \left( H_0^1(\Omega) \right)^2 \setminus \{0 \}} }  \frac{ \displaystyle{\int_{\Omega} \left[  \nabla a(x) \cdot (x-x_0) {\vert \nabla u(x) \vert}^2 + \nabla b(x) \cdot (x-x_0){\vert \nabla v(x) \vert}^2 \right] dx}}{\displaystyle{{  \int_{\Omega}} u(x) v(x) dx } } $$
which is a contradiction.

Then the problem \eqref{BN} does not admit solutions. Consequently we obtain the desired result.   \hfill $\blacksquare$

\subsection{Estimates of $\omega(a, b)$}
We end this section, using the techniques in \cite{HY},  we give some estimates of $\omega(a, b). $
\begin{proposition} \mbox{}
\begin{itemize}
\item[(1)] We assume that $ a, b \in C^1(\Omega) $ and there exists $z_0 \in \Omega $ such that $ \tilde{a}(z_0)  + \tilde{b}(z_0) < 0, $ then $ \omega(a, b) = - \infty. $
\item[(2)] We assume that $ a$ respectively $b \in H^1(\Omega) \cap C(\bar{\Omega}) $ satisfying $\eqref{1.3}$ respectively  $\eqref{1.4}$ and  we have
\item[(2.i)] If $ k > 2, l > 2$ and $ a, b \in C^1(\Omega), $ then $ \omega(a, b) = 0. $
\item[(2.ii)] If $ k = 2, l > 2 $ or $ k > 2, l = 2, $ and  $a$, $b$  satisfying moreover $\eqref{1.5}$ and $\eqref{1.6}$ with $ \tilde{a}(x) \geq 0$, $\tilde{b}(x) \geq 0 $ a.e. $ x \in \Omega, $ then
\begin{equation} \label{Inq1}
\frac{N^2}{16} \min \left( A_2 ,  l B_l ( \diam \ \Omega )^{l-2} \right)  \leq \omega(a, b) \leq \frac{A_2}{2} \lambda_1 (\diam \ \Omega )^2
\end{equation}
 and
\begin{equation} \label{Inq2}
 \frac{N^2}{16} \min \left( k A_k (\diam \ \Omega )^{k-2} , B_2 \right) \leq \omega(a, b) \leq \frac{B_2}{2} \lambda_1 (\diam \ \Omega )^2
 \end{equation}
 \item[(2.iii)] If $ 0 < k \leq 2, 0 < l \leq 2, $ $ a, b $ satisfy the conditions $\eqref{1.5}$ and $\eqref{1.6}$ respectively, then $$ \frac{N^2}{16} \min \left( k A_k ( \diam \ \Omega )^{k-2},  l B_l (\diam \ \Omega )^{l-2} \right) \leq \omega(a, b). $$
\end{itemize}
\end{proposition}
\
\textbf{Proof}
Let $ \varphi \in C_c({\mathbb{R}}^N)$ such that $0 \leq \varphi \leq 1 $ on $ {\mathbb{R}}^N, $ $ \varphi \equiv 1 $ on $ B(0,2r) $ and $ \varphi \equiv 0 $ on $ {\mathbb{R}}^N \setminus B(0,2r), $ where $0< r < 1. $
Set $ \varphi_j(x)= \varphi(j(x-z_0))$ for $ j \in {\mathbb{N}}^{*}, $ we have
$$
\begin{array}{cll}
\omega(a, b) &\leq&  \displaystyle{ \frac{1}{4} } \frac{ \displaystyle{ \int_{\Omega} \left( \tilde{a}(x) + \tilde{b}(x) \right) \vert \nabla {\varphi_j(x) \vert}^2 dx}}{ \displaystyle{ \int_{\Omega} \varphi_j^2(x)  dx }} \\
&\leq&  \displaystyle{ \frac{1}{4}}  \frac{  \displaystyle{ \int_{B(b,\frac{2r}{j})} \left( \tilde{a}(x) + \tilde{b}(x) \right) \vert \nabla {\varphi_j(x) \vert}^2 dx}}{  \displaystyle{\int_{B(b,\frac{2r}{j})} \varphi_j^2(x)  dx }.} \\
\end{array}
$$
Using the change of variable $y=j(x-z_0), $ we get
$$
\omega(a, b) \leq \displaystyle{ \frac{j^2}{4} } \frac{ \displaystyle{ \int_{B(0,2r)} \left( \tilde{a}(\frac{y}{j}+ z_0) + \tilde{b}(\frac{y}{j}+ z_0) \right) \vert \nabla {\varphi(y) \vert}^2 dy}}{ \displaystyle{ \int_{B(0,2r)} \varphi^2(y) dy }}.
$$
Applying the Dominated Convergence Theorem, we obtain the desired result when $j$ goes to infinity. \\
Now we will prove $(2.i). $ \\
Using \eqref{1.3}, \eqref{1.4} and since $a, b \in C^1(\Omega)$ in a neighborhood $ V $ of $x_1, $ we write
\begin{equation} \label{3.10}
 a(x) = a(x_1) + {A_k\vert x-x_1 \vert}^k + \theta_{a}(x)
\end{equation}
\begin{equation}  \label{3.11}
 b(x) = b(x_1) + {B_l \vert x-x_1 \vert}^l + \theta_{b}(x),
\end{equation}
where $\theta_{a}(x)$ and $\theta_{b}(x) \in C^1(V)$ are such that
\begin{equation} \label{3.12}
 \lim_{x \rightarrow x_1} \frac{\vert \theta_{a}(x) \vert}{{\vert x-x_1 \vert}^k } =0 \quad \mbox{and} \quad \lim_{x \rightarrow x_2} \frac{\vert \theta_{b}(x) \vert}{{\vert x-x_1 \vert}^l } =0.
 \end{equation}
From \eqref{3.12}, we get the existence of $r,$ such that $0< r < 1$ and
\begin{equation} \label{3.13}
 \vert \theta_{a}(x) \vert \leq {\vert x-x_0 \vert}^k \quad \mbox{and}
  \quad \vert \theta_{b}(x) \vert \leq {\vert x-x_0 \vert}^l,  \mbox{for all} \quad  x \in B(x_0,2r)  \subset V.
 \end{equation}
Let $ \varphi_{j}(x)= \varphi(j(x-x_0))$ define as in the proof of $(1); $ we have
$$
0 \leq \omega(a, b) \leq \frac{1}{4}  \frac{ \displaystyle{ \int_{\Omega} \left( \tilde{a}(x) +  \tilde{b}(x) \right) \vert \nabla \varphi_{j}(x)\vert^2 dx}}{ \displaystyle{\int_{\Omega} \varphi_{j}^2(x) dx}}.
$$
Using \eqref{3.10} and \eqref{3.11}, we obtain
$$
\begin{array}{cll}
0 \leq \omega(a, b) &\leq&  \displaystyle{ \frac{1}{4} } \frac{\displaystyle{  \int_{B(x_0,\frac{2r}{j})} \left(  k A_k {\vert x-x_0 \vert}^k + l B_l {\vert x-x_0 \vert}^l \right) \vert \nabla \varphi_{j}(x)\vert^2 dx }}{ \displaystyle{\int_{B(x_0,\frac{2r}{j})} \varphi_{j}^2(x) dx } } \\ &+&  \displaystyle{  \frac{1}{4} } \frac{\displaystyle{ \int_{B(x_0,\frac{2r}{j})} \left( \nabla \theta_{a}(x)\cdot(x-x_0 ) + \nabla \theta_{b}(x)\cdot ( x-x_0 ) \right)  \vert \nabla \varphi_{j}(x)\vert^2 dx}}{ \displaystyle{ \int_{B(x_0,\frac{2r}{j})} \varphi_{j}^2(x)  dx }. }   \\
\end{array}
$$
By a simple change of variable $ y=j(x-x_0)$ and integrating by parts, we obtain
$$
\begin{array}{cll}
0 \leq \omega(a, b) &\leq  & \displaystyle{ \frac{k A_k}{4 j^{k-2}} } 
\frac{\displaystyle{\int_{B(0,2r)}{\vert y \vert}^k \vert \nabla \varphi(y)\vert^2 dy}}{\displaystyle{ \int_{B(0,2r)} \varphi(y)^2  dy}}
- \displaystyle{ \frac{j}{4}} \frac{\displaystyle{  \int_{B(0,2r)} \theta_a\Big(\frac{y}{j}+x_0\Big)\cdot\div(y \vert \nabla \varphi(y)\vert^2) dy}}{\displaystyle{ \int_{B(0,2r)} {\varphi(y)}^2 dy}} + \\
&+& \displaystyle{ \frac{l B_l}{4 j^{l-2}} }\frac{ \displaystyle{ \int_{B(0,2r)} {\vert y \vert}^l \vert \nabla \varphi(y)\vert^2 dy }}{ \displaystyle{ \int_{B(0,2r)} {\varphi(y)}^2 dy}} - \displaystyle{  \frac{j}{4}} 
\frac{ \displaystyle{ \int_{B(0,2r)} \theta_b(\frac{y}{j}+x_0)\cdot\div \left( y {\vert \nabla \varphi(y)\vert}^2 \right) dy}}{\displaystyle{ \int_{B(0,2r)} {\varphi(y)}^2 dy}}.   \\
\end{array}
$$
Using \eqref{3.13}, we get
$$
\begin{array}{cll}
0 \leq \omega(a, b) &\leq & \displaystyle{ \frac{k A_k}{4 j^{k-2}} } \frac{ \displaystyle{\int_{B(0,2r)}{\vert y \vert}^k \nabla \varphi(y) dy}}{\displaystyle{ \int_{B(0,2r)}  \varphi^{2}(y) dy} } + \displaystyle{  \frac{C}{j^{k-1}}} \frac{\displaystyle{ \int_{B(0,2r)} {\vert y \vert }^k dy}}{\displaystyle{\int_{B(0,2r)} \varphi^{2}(y) dy}} \\
&+& \displaystyle{ \frac{l B_l}{4 j^{{l-2}}}} \frac{ \displaystyle{ \int_{B(0,2r)} \vert y \vert^l  \vert \nabla \varphi(y) \vert^2 dy }}{\displaystyle{ \int_{B(0,2r)} \varphi^{2}(y) dy}} +  \displaystyle{ \frac{C}{j^{l-1}}} \frac{ \displaystyle{ \int_{B(0,2r)} \vert y \vert^l dy}}{\displaystyle{\int_{B(0,2r)} \varphi^{2}(y)  dy}},
\end{array}
$$
where $ C= \displaystyle{ \max_{y \in B(0,2r)} \vert \div(y \vert \nabla \varphi(y) \vert^2 )\vert.} $

Therefore, for $k > 2$ and $l > 2$,  when $j$  tends to $\infty$ we obtain  $\omega_(p, q) = 0$.

{ To prove $(2.ii), $ we start by the case $k = 2$ and $l > 2. $ We show the left-hand inequality in \eqref{Inq1}.
Since $a$ and $b$ satisfy $\eqref{1.5}$ and $\eqref{1.6}$, respectively, for all pairs $(u, v) \in E \setminus \{0 \}, $ we have
\begin{align*}
\phi_{a,b}(u,v) &\geq \frac{1}{2} A_2 \frac{\displaystyle{\int_{\Omega} \vert ( x-x_0 ) \cdot \nabla u(x) \vert^2 dx}}{\displaystyle{\int_{\Omega} u(x)v(x) dx}} + \frac{l}{4} B_l \frac{\displaystyle{\int_{\Omega} \vert x-x_0\vert^{l-2} \vert ( x-x_0 ) \cdot \nabla v(x) \vert^2 dx}}{\displaystyle{\int_{\Omega} u(x)v(x) dx}} \\
&\geq  \frac{1}{2} A_2  \frac{\displaystyle{\int_{\Omega} \vert ( x-x_0 ) \cdot \nabla u(x)\vert^2 dx}}{\displaystyle{\int_{\Omega} u(x)v(x) dx}} + \frac{l}{4} B_l (\diam \ \Omega )^{l-2} \frac{\displaystyle{ \int_{\Omega} 
 \vert ( x-x_0 ) \cdot \nabla v(x) \vert^2 dx}}{\displaystyle{\int_{\Omega} u(x)v(x) dx}}.
\end{align*}
}
By applying Lemma 2.1 in {\cite{HY}} for $ t=0, $ we find
\begin{align}
\phi_{a,b}(u,v) \geq & \frac{1}{2} A_2  \left( \frac{N}{2} \right)^2 \frac{\displaystyle{\int_{\Omega}  \vert u(x) \vert^2 dx}}{\displaystyle{\int_{\Omega} u(x)v(x)dx}}+ \frac{l}{4} B_l (diam \ \Omega )^{l-2} \left( \frac{N}{2} \right)^2 \frac{\displaystyle{\int_{\Omega}  \vert v(x) \vert^2 dx}}{\displaystyle{\int_{\Omega} u(x)v(x) dx}} \nonumber.
\end{align}
This implies that
$$ \displaystyle{\omega(a, b) \geq \frac{N^2}{16} \min( A_2 , l B_l (\diam \ \Omega )^{l-2}). }  $$
Similarly, we deduce in the case $ k > 2 $ and $ l=2, $ that
 $$ \displaystyle{\omega(a, b) \geq \frac{N^2}{16} \min( k A_k (\diam \ \Omega )^{k-2} , B_2  ). }  $$
Now we prove the right-hand inequality in \eqref{Inq1} and \eqref{Inq2}. Let $ \psi_{j}(x)= \varphi_1(j(x-x_0))$ for $j \in \mathbb{N}$ large enough, where $\varphi_1$ is a positive
eigenfunction corresponding to the first eigenvalue $\lambda_1$ of the operator $-\Delta$ in $  H_0^1(\Omega). $ \\
We have
$$
0 \leq \omega(a, b) \leq \frac{1}{4}  \frac{ \displaystyle{ \int_{\Omega} (\tilde{a}(x) + \tilde{b}(x) ) \vert \psi_j(x)\vert^2 dx}}{\displaystyle{ \int_{\Omega} \psi_{j}^2(x)  dx }}.
$$
 Using $(4)$ and $(5)$, we obtain
$$
\begin{array}{cll}
0 \leq \omega(a, b) &\leq&   \frac{\displaystyle{ \int_{x_0+\frac{1}{j} \Omega} \left(  2 A_2 {\vert x-x_0 \vert}^2 + l B_l {\vert x-x_0 \vert}^l \right) \vert \nabla \psi_{j}(x)\vert^2 dx}}{ \displaystyle{4\int_{x_0,\frac{1}{j} \Omega} \psi_{j}^2(x) dx } } \\ &+&     \frac{ \displaystyle{ \int_{x_0+\frac{1}{j} \Omega} \left( \nabla \theta_{a}(x)\cdot(x-x_0 ) +  \nabla \theta_{b}(x)\cdot ( x-x_0 )\right) \vert \nabla \psi_{j}(x) \vert^2 dx}}{ \displaystyle{4\int_{x_0+\frac{1}{j} \Omega} \psi_{j}^2(x)  dx }.  } \\
\end{array}
$$
By a simple change of variable $ y=j(x-x_0)$ and integrating by parts, we have by~\eqref{3.12} 

$$
\begin{array}{cll}
0 \leq \omega(a, b) &\leq  &  \displaystyle{ \frac{A_2}{2} } \frac{ \displaystyle{ \int_{\Omega}{\vert y \vert}^2 \vert \nabla \varphi_1(y)\vert^2 dy}}{\displaystyle{ \int_{\Omega} \varphi_1(y)^2  dy}} + \displaystyle{ \frac{C}{j} }
\frac{ \displaystyle{\int_{\Omega} \vert y \vert^k  dy}}{\displaystyle{\int_{\Omega} {\varphi_1(y)}^2 dy}},\\
&+&  \displaystyle{\frac{l B_l}{4 j^{l-2}}} \frac{ \displaystyle{ \int_{\Omega} {\vert y \vert}^l \vert \nabla \varphi_1(y)\vert^2 dy }}{ \displaystyle{\int_{\Omega} {\varphi_1(y)}^2 dy}}
+   \frac{C}{j^{l-1}} \frac{ \displaystyle{\int_{\Omega} \vert y \vert^l  dy}}{\displaystyle{\int_{\Omega} {\varphi_1(y)}^2 dy}},   \\
\end{array}
$$
where 
$ C = {\displaystyle{ \max_{y \in \Omega}}}\  \vert\div(y \nabla \varphi_1(y) \vert^2) \vert. $ Letting $ j \rightarrow \infty $ we get
$$ 0 \leq \omega(a, b) \leq  \frac{ A_2}{2}  \frac{ \displaystyle{ \int_{\Omega}{\vert y \vert}^2 \vert \nabla \varphi_1(y)\vert^2 dy}}{\displaystyle{ \int_{\Omega} \varphi_1(y)^2  dy}} , $$
therfore
$$ 0 \leq \omega(a, b) \leq \frac{ A_2}{2} \lambda_1 (\diam \ \Omega )^2. $$
Similarly, we deduce in the case $ k > 2 $ and $ l=2, $ that
 $$ 0 \leq \omega(a, b) \leq \frac{ B_2}{2} \lambda_1 (\diam \ \Omega )^2. $$
\\
\color{black}{Let us now prove $(2.iii). $ Since $a$ and $b$ satisfy $\eqref{1.5}$ and $\eqref{1.6}$, respectively,\\ for all $(u, v) \in E \setminus \{0 \}, $ we have
\begin{align*}
\phi_{a,b}(u,v) &\geq \frac{k}{4} A_k \frac{\displaystyle{\int_{\Omega} \vert x-x_0\vert^{k-2} \vert ( x-x_0 ) \cdot \nabla u \vert^2 dx}}{\displaystyle{\int_{\Omega} u v dx}} + \frac{l}{4} B_l \frac{\displaystyle{\int_{\Omega} \vert x-x_0\vert^{l-2} \vert ( x-x_0 ) \cdot \nabla v \vert^2 dx}}{\displaystyle{\int_{\Omega} u v dx}} \\
&\geq  \frac{k}{4} A_k ( \diam \ \Omega )^{k-2} \frac{\displaystyle{\int_{\Omega} \vert ( x-x_0 ) \cdot \nabla u \vert^2 dx}}{\displaystyle{\int_{\Omega} u v dx}} + \frac{l}{4} B_l (\diam \ \Omega )^{l-2} \frac{\displaystyle{ \int_{\Omega}  \vert ( x-x_0 ) \cdot \nabla v \vert^2 dx}}{\displaystyle{\int_{\Omega} u v dx}}.
\end{align*}
}
By applying Lemma 2.1 in {\cite{HY}} for $ t=0, $ we find
\begin{align}
\phi_{a,b}(u,v) \geq & \frac{k}{4} A_k (\diam \ \Omega )^{k-2} \left( \frac{N}{2} \right)^2 \frac{\displaystyle{\int_{\Omega}  \vert u \vert^2 dx}}{\displaystyle{\int_{\Omega} u v dx}}+ \frac{l}{4} B_l (\diam \ \Omega )^{l-2} \left( \frac{N}{2} \right)^2 \frac{\displaystyle{\int_{\Omega}  \vert v \vert^2 dx}}{\displaystyle{\int_{\Omega} u v dx}} \nonumber.
\end{align}
This implies that
$$ \displaystyle{\omega(a, b) \geq \frac{N^2}{16} \min(k A_k (\diam \ \Omega )^{k-2}, l B_l (\diam \ \Omega )^{l-2})}. $$ \hfill $\blacksquare$
   \bibliographystyle{plain}

 \end{document}